\documentclass[submission]{dmtcs}

\usepackage[latin1]{inputenc}
\usepackage{subfigure}
\usepackage[round]{natbib}

\RCSdef$Revision: 1.3 $\endRCSdef
\rcsMajMin

\def\b{\beta} 

\def\G{\Gamma}
\def\D{\triangle}

\def\d{\delta}
\def\o{\omega}  
 
\def\l{\lambda}

\def\m{\mu}

\def\v{\vdash}
 
\def\p{\succ}
 
\def\<{\langle}
\def\>{\rangle}
\def\F{\displaystyle\frac}
\def\red{\triangleright^*}
\def\rd{\triangleright}
\def\ss{\sqsupset}
\def\sq{\sqsupseteq}

\newtheorem{theo}{Theorem}[section]
\newtheorem{lm}{Lemma}[section]
\newtheorem{re}{Remark}[section]

\newtheorem{df}{Definition}[section]

\newtheorem{ex}{Example}[section]

\author{Karim Nour\addressmark{1}\thanks{knour@univ-savoie.fr} \and
  Khelifa Saber\addressmark{2} \thanks{ksabe@univ-savoie.fr}}

\title[The call-by-value $\l \m ^{\wedge \vee}$-calculus]{Confluency property of the call-by-value $\l \m ^{\wedge \vee}$-calculus}

\address{\addressmark{1} LAMA - Equipe de logique , Universit\'e de
Savoie ,  F-73376 Le Bourget du Lac, France\\
\addressmark{2} LAMA - Equipe de logique , Universit\'e de
Savoie ,  F-73376 Le Bourget du Lac, France}
\keywords{Call-by-value, Church-Rosser, Propositional classical logic,
Parallel reduction, Complete development}

\revision{\rcsMaj}

\received{}
\revised{}
\accepted{}

\begin{document}

\maketitle
\begin{abstract} 
In this paper, we introduce the $\l \m ^{\wedge \vee}$- call-by-value
calculus and we give a proof of the Church-Rosser property of this
system. This proof is an adaptation of that of \cite{and} which uses an extended parallel reduction method and
complete development.
\end{abstract}

\clearpage
\tableofcontents

\section{Introduction}

 $\;\;\;$\cite{gen} introduced the natural deduction
 system  to study the notion of proof. The full classical natural
 deduction system is well adapted for the human reasoning. By full we
 mean that all the connectives ($\to$, $\wedge$ and $\vee$) and $\bot$ (for the absurdity) are
 considered as primitive. As usual, the negation is defined by $\neg A = A \to
 \bot$. Considering this logic from the computer science of view is
 interesting because, by the Curry-Howard correspondence, formulas can
 be seen as types for the functional programming languages and correct
 programs can be extracted. The corresponding
 calculus is an extension of M. Parigot's $\l \m$-calculus with
 product and coproduct, which is denoted by $\l \m ^{\wedge
 \vee}$-calculus.

 \cite{deG2} introduced the  typed $\l \m ^{\wedge
\vee}$-calculus to code the classical natural deduction system, and
showed that it enjoys the main important properties: the strong normalization, the confluence and
the subformula property. This would guarantee that proof
normalization  may be interpreted as an evaluation process. As far as
we know the typed $\l \m ^{\wedge \vee}$-calculus is the first extension of
the simply typed $\l$-calculus which enjoys all the above
properties. \cite{Rit1} introduced  an extension of the $\l \m$-calculus that
features disjunction as primitive (see also \cite{Rit2}). But their
system is rather different since they take as primitive a classical
form of disjunction that amounts to $\neg A \to B$. Nevertheless,
\cite{Rit3} give another extension of the $\l
\m$-calculus with an intuitionistic disjunction. However, the
reduction rules considered are not sufficient to guarantee that the
normal forms satisfy the subformula property. The question of the
strong normalization of the full logic has interested several authors,
thus one finds in \cite{davnou1}, \cite{matt} and \cite{Nour7} different proofs of this result.

 From a computer science point of view, the $\l \m ^{\wedge
 \vee}$-calculus may be seen as the kernel of a typed call-by-name
 functional language featuring product, coproduct and control
 operators. However we cannot apply an arbitrary reduction for
 implementation of programming languages, we have to fix a reduction
 strategy and usually it is the call-by-value strategy.  Many programming langagues and control operations were
 developed  through the studies of the call-by-value
 variant like {\tt ML} and {\tt Lisp} for $\l$-calculus, the calculus of exception
 handling {\tt $\l_{exn}^{\to}$} and {\tt $\m$PCF$_V$} for the $\l
 \m$-calculus.
 \cite{Ong} showed that {\tt $\m$PCF$_V$} is sufficiently strong to express the
 various control constructs such as the {\tt ML-style raise} and the
 first-class continuations {\tt callcc}, {\tt throw} and {\tt abort}. In this
 sense, it seems to be important to study the call-by-value version of
 $\l \m ^{\wedge \vee}$-calculus.

Among the important properties required in any abstract reduction
system, there is the confluence which ensures the uniqueness of the
normal form (if it exists). The notion of parallel reduction which is based
on the method of Tait and Martin-Löf is a good tool to prove the
confluence property for several reduction systems. The idea is very
clear and intuitive: It consists in reducing
a number of redexes existing in the term simultaneously. However,
this method does not work for the $\l \m ^{\wedge \vee}$-calculus. In
fact the diamond property which stipulates that: If $t \succ t'$ then $ t'
\succ t^* $ (where $t^*$ is usually referred as the complete
development of $t$) does not hold because more complicated situations
appear, and that is due to  the presence of the permutative reductions
``$((u\,[x.v,y.w])\,\varepsilon) \rd (u\,[x.(v\,\varepsilon),y.(w\,\varepsilon)])$''. Hence the proof of the confluence
becomes hard and not at all trivial as it seems to be.

 Consider the terms $t=(((u\,[x.v,y.w])\,[r.p,s.q])\,\varepsilon)$, $t_1= ((u\,[x.(v\,[r.p,s.q]),y.(w\,[r.p,s.q])])\,\varepsilon)$,
 and  $t_2= ((u\,[x.v,y.w])\,[r.(p\,\varepsilon),s.(q\,\varepsilon)])$. We have: $t \succ t_1$ and $t \succ t_2$,
 if we want the diamond property to hold, $t_1$ and $t_2$ must be
 reduced to the same term $t^*$ by one reduction step , however this
 is not possible. To make it possible we need another step of
 permutative reduction. We consider such a successive sequence of
 reductions as a one parallel reduction step, i.e, we
 follow the permutative reductions in the term step by step to a certain depth
 which allows  to join and consider this sequence as a one reduction step. The notion of Prawitz' s {\em segment} yields the formulation of  this new parallel reduction. Therefore the difficulties are overcome by extending this notion to our
system (see \cite{andd}, \cite{and}, \cite{Pra1} and
\cite{Pra2}) and considering the extended structural reductions along
this {\it segment} which allow us to define a complete development to
obtain directly the common reductum, hence the Church-Rosser
property. This is exactly what is done in \cite{and}; our proof is
just a checking that this method is well adapted to provide the
diamond property for the call-by-value $\l \m ^{\wedge \vee}$-calculus
including the  symmetrical rules. Thus $t_1 \succ t^*$
and $t_2 \succ t^*$, where $t^* =(u\,[x.(v\,[r.(p\,\varepsilon),s.(q\,\varepsilon)]),y.(w\,[r.(p\,\varepsilon),s.(q\,\varepsilon)])])$.

The paper is organized as follows. Section 2 is an introduction to
 the typed system, the relative cut-elimination procedure of $\l \m
 ^{\wedge \vee}$-calculus and the call-by-value $\l \m ^{\wedge
 \vee}$-calculus. In section 3, we define the parallel reduction
 related to the notion of {\it segment-tree}, thus we give the key lemma
 from which the diamond-property is directly deduced. Section 4 is
 devoted to the proof of the key lemma. We conclude with some future work.

\section{Notations and  definitions}
\begin{df}
We use notations inspired by \cite{and}. 
\begin{enumerate}
\item Let $\mathcal{X}$ and $\mathcal{A}$ be two disjoint alphabets
for distinguishing the $\lambda$-variables and $\mu$-variables
respectively.  We code deductions by using a set of terms
$\mathcal{T}$ which extends the $\l$-terms and is given by the
following grammar (which gives terms at the untyped level):
\begin{center}
$
\mathcal{T} \; := \;\mathcal{X} \; | \; 
\lambda\mathcal{X}.\mathcal{T}\;
|\; (\mathcal{T}\;\;\mathcal{E}) \; | \; \< \mathcal{T},\mathcal{T} \> 
\; | \;$$
\omega_1 \mathcal{T}$$ \; |\;$$ \omega_2 \mathcal{T}$$ \; | \;
\mu\mathcal{A}.\mathcal{T} \; | \; (\mathcal{A}\; \; \mathcal{T})
$

$
\mathcal{E} \; := \; \mathcal{T} \; | \; $$\pi_1$$ \; | \;$ $\pi_2 $$\; 
|
\; [\mathcal{X}.\mathcal{T},\mathcal{X}.\mathcal{T}]
$
\end{center}

An element of the set $\mathcal{E}$ is said to be an 
$\mathcal{E}$-term. Application between two $\mathcal{E}$-terms $u$
and $\varepsilon$ is denoted by $(u \;\varepsilon)$.

\item The meaning of the new constructors is given by the typing rules
below where $\G$ (resp. $\Delta$) is a context, i.e. a set of declarations of the form
$x : A$ (resp. $a :  A$) where $x$ is a $\l$-variable (resp. $a$
is a $\m$-variable) and $A$ is a formula.

\begin{center}
 $\F{}{\Gamma, x:A\,\, \vdash x:A\,\, ; \, \Delta}{ax}$
\end{center}

\begin{center}
$\F{\Gamma, x:A \vdash t:B;\Delta}{\Gamma \vdash \lambda x.t:A \to 
B;\Delta}{\to_i}
\quad\quad\quad 
\F{\Gamma \vdash u:A \to B;\Delta \quad \Gamma \vdash
v:A;\Delta}{\Gamma\vdash (u\; \; v):B;\Delta}{\to_e}$
\end{center}

\begin{center}
$\F{\Gamma \vdash u:A;\Delta \quad \Gamma \vdash v:B ; \Delta}{\Gamma
\vdash \<u,v\>:A \wedge B ; \Delta}{\wedge_i}$
\end{center}

\begin{center}
$\F{\Gamma \vdash t:A \wedge B ; \Delta}{\Gamma \vdash (t\;\;\pi_1):A ;
\Delta}{\wedge^1_e} 
\quad 
\F{\Gamma \vdash t:A\wedge B ; \Delta}{\Gamma \vdash (t\;\;\pi_2):B ; 
\Delta}{\wedge^2_e}$
\end{center}

\begin{center}
$\F{\Gamma\vdash t:A;\Delta}{\Gamma\vdash \omega_1
    t:A \vee B ;\Delta}{\vee^1_i}
\quad
\F{\Gamma \vdash t:B; \Delta}{\Gamma\vdash \omega_2 t:A\vee B 
;\Delta}{\vee^2_i}$
\end{center}

\begin{center}
$\F{\Gamma \vdash t:A\vee B ;\Delta\quad\Gamma, x:A \vdash u:C ;
\Delta\quad\Gamma, y:B \vdash v:C ; \Delta}{\Gamma \vdash (t\;\;[x.u, 
y.v]):C ; \Delta}{\vee_e}$ 
\end{center}

\begin{center}
$\F{\Gamma\vdash t:A ;\Delta, a:A}{\Gamma \vdash (a\;\;t):\bot ;\Delta,
a:A}{\bot_i}
\quad
\F{\Gamma\vdash t:\bot; \Delta, a:A}{\Gamma \vdash\mu a.t:A; 
\Delta}{\bot_e}$
\end{center}

\item A term in the form $(t\;[x.u,y.v])$ (resp $\mu a.t$) is called
  an ${\vee_e}$-term (resp ${\bot_e}$-term).

\item The cut-elimination procedure corresponds to the  reduction rules 
given below. They are those we need to the subformula 
property.

\begin{itemize}
\item $(\lambda x.u \;\; v) \triangleright_{\beta} u[x:=v]$

\item $(\<t_1,t_2\>\;\;\pi_i) \triangleright_{\pi} t_i$

\item $(\omega_i t\;\;[x_1.u_1,x_2.u_2]) \triangleright_{D}  u_i[x_i:=t]$

\item $((t\;\;[x_1.u_1,x_2.u_2])\;\;\varepsilon) \triangleright_{\delta}
(t\;\;[x_1.(u_1\; \varepsilon),x_2.(u_2\;\varepsilon)])$

\item $(\m a.t\;\; \varepsilon) \triangleright_{\mu} \m
  a.t[a:=^*\varepsilon]$

where $ t[a:=^*\varepsilon]$ is obtained from $t$ by replacing
inductively each subterm in the form $(a \; v)$ by $(a \; (v \; 
\varepsilon))$.
\end{itemize}

\item Let $t$ and $t'$ be terms. The notation $t
\triangleright t'$ means that $t$ reduces to $t'$ by using one step
of the reduction rules given above. Similarly, $t \triangleright^* t'$
means that $t$ reduces to $t'$ by using some steps of the reduction
rules given above.
\end{enumerate}
\end{df}

The following result is straightforward

\begin{theo}(Subject reduction)
If $\G \vdash t : A ; \Delta$ and $t \triangleright^* t'$, then $\G 
\vdash t' : A ; \Delta$.
\end{theo}

We have also the following properties (see \cite{and}, \cite{davnou1},
 \cite{deG2}, \cite{matt}, \cite{Nour7} and \cite{nour}).
\begin{theo}(Confluence) If $t\triangleright^* t_1$ and $t\red t_2$, 
then there exists $t_3$ such that $t_1\red t_3$ and $t_2\red t_3$.
\end{theo}

\begin{theo}(Strong normalization) If $\G \vdash t : A ; \Delta$, 
then $t$ is strongly normalizable.
\end{theo}

\begin{re}
Following the call-by-value evaluation discipline, in an
application the evaluator has to diverge if the argument diverges.
For example, in the call-by-value $\l$-calculus, we are allowed to reduce the
$\b$-redex $(\l x. u \; v)$ only when $v$ is a value. In
$\l\mu$-calculus, the terms $\mu a.u$ and $(u\;[x_1.u_1,x_2.u_2])$
cannot be taken as values, then the terms $(\l x. t \; \mu a.u)$ and
$(\l x.t \; (u\;[x_1.u_1,x_2.u_2]))$ cannot be reduced. This will be
able to prevent us from reaching many normal forms. To solve this
problem, we introduce symmetrical rules (${\d_v'}$ and ${\mu_v'}$)
allowing to reduce these kinds of redexes.
\end{re}

Now we introduce the call-by-value version of the $\l \m ^{\wedge \vee
}$-calculus. From a logical point of view a value corresponds to
an introduction of a connective; this is the reason why the Parigot's
naming rule is considered as the introduction rule of $\bot$.

\begin{df}
\begin{enumerate}
\item The set of values $\mathcal{V}$ is given by the following grammar:
\begin{center}
$
\mathcal{V} \; := \;\mathcal{X} \; | \; \lambda\mathcal{X}.\mathcal{T}\;
| \; \< \mathcal{V},\mathcal{V} \> \; | \;$$
\omega_1 \mathcal{V}$$ \; |\;$$ \omega_2 \mathcal{V}$$ \; | \; (\mathcal{A}\; \; \mathcal{T})
$
\end{center}
Values are denoted $U,V,W,...$

\item The reduction rules of the call-by-value $\l
\m^{\wedge \vee}$-calculus  are the followings:

\begin{itemize}

\item $(\lambda x.t \;\; V) \triangleright_{\b_v} t[x:=V]$

\item $(\<V_1,V_2\>\;\;\pi_i) \triangleright_{\pi_v} V_i$

\item $(\omega_i V\;\;[x_1.t_1,x_2.t_2]) \triangleright_{D_v}
t_i[x_i:=V]$

\item $((t\;\;[x_1.t_1,x_2.t_2])\;\;\varepsilon) \triangleright_\d
(t\;\;[x_1.(t_1\; \varepsilon),x_2.(t_2\;\varepsilon)])$

\item $(V\;(t\;\;[x_1.t_1,x_2.t_2])) \triangleright_{\d_v'}
(t\;\;[x_1.(V\;t_1),x_2.(V\;t_2)])$

\item $(\m a.t\;\; \varepsilon) \triangleright_{\m} \m
a.t[a:=^*\varepsilon]$

\item $(V\;\;\m a.t) \triangleright_{\m'_v} \m a.t[a:=_*V]$

where $t[a:=_* V]$ is obtained from $t$ by replacing
 inductively each subterm in $t$ in the form $(a\;u)$ by
 $(a\;(V\;u))$.
\end{itemize}

The first three rules are called logical rules and the others are
 called structural rules.

\item The one-step reduction $\rd_v$ of the call-by-value $\l \m ^{\wedge\vee}$-calculus is defined as the union of the seven rules given
 above. As usual $\red_v$ denotes the transitive and reflexive closure of
 $\rd_v$. 
\end{enumerate}
\end{df}

 The following lemma expresses the fact that the set of values is
 closed under reductions. In the remainder of this paper, this fact
 will be used implicitly.

\begin{lm} 
If $V$ is a value and $V \red_v W$, then $W$ is a value.
\end{lm}

\begin{proof}
From the definition of the set of values.
\end{proof}

\begin{theo}(Subject reduction)
 If $\G \v t:A\; ;\D$ and $t \red_v t'$, then $\G \v t':A\; ;\D$.
\end{theo}

\begin{proof} Since the reduction rules correspond to the
  cut-elimination procedure, we check easily that the type is preserved from the redex to
  its reductom.\end{proof}

The rest of this paper is an extention of \cite{and} to our calculus
according to the new considered reduction rules $\d'_v$ and $\m'_v$. One
can find all the notions given here in \cite{and}. Since the new
symmetrical rules that we add don't create any critical pair with the
existing rules, then in the examples and proofs that we give, one
will mention only the cases related to these new rules and  check that they
don't affect  the core  of \cite {and}'s work.

\section{The extended structural reduction}

\begin{df}
\begin{enumerate}
\item Let $t$ be a term, we define a binary relation denoted by
  $\ss_t$ on subterms of $t$ as follows:
\begin{itemize}

\item $(u\;\;[x_1.u_1,x_2.u_2]) \ss_t u_i$

\item $\m a.u \ss_t v$, where $v$ occurs in $u$ in the form $(a\;v)$

\end{itemize} 

If $u \ss_t v$ holds, then $v$ is called a segment-successor of $u$,
and  $u$ is called a segment-predecessor of $v$. We denote by $\sq_t$
the reflexive and transitive closure of $\ss_t$.

\item Let $r$ be a subterm of a term $t$, such that $r$ is a $\vee_e$- or
  $\bot_e$-term and $r$ has no segment-predecessor in $t$. A segment-tree
  from $r$ in $t$ is a set ${\cal O}$ of subterms of $t$, such
  that for each $w \in {\cal O}$:
\begin{itemize}
\item $r \sq_t w$
\item $w$ is a $\vee_e$- or $\bot_e$-term
\item For each subterm $s$ of $t$, such  that $r \sq_t s\sq_t w$ then
  $s \in {\cal O}$
\end{itemize}
$r$ is called the root of ${\cal O}$.
\item Let ${\cal O}$ be a segment-tree from $r$ in $t$, a subterm $v$ of
  $t$ is called an acceptor  of ${\cal O}$ iff $v$ is a
  segment-successor of an element of ${\cal O}$ and $v $ is not in ${\cal
  O}$.

\item A segment-tree ${\cal O}$ from $r$ in $t$ is called the maximal
  segment-tree iff no acceptor of ${\cal O}$ has a segment successor
  in $t$.

\item The acceptors  of ${\cal O}$ are indexed by the letter ${\cal O}$.

\item Let ${\cal O}$ be a segment-tree from $t$ in $t$ itself, and $t \red_v
  t'$, then we define canonically a corresponding segment-tree to
  ${\cal O}$ in $t'$ by the transformation of indexes from redexes to
  their residuals. This new segment-tree is denoted also by ${\cal
  O}$ if there is no ambiguity.

\end{enumerate}
\end{df}

\begin{re} 
For typed terms, all the elements of a segment-tree have the
same type.
\end{re}

\begin{df}
 Let ${\cal O}$ be a segment-tree from $r$ in $t$, suppose that $r$
 occurs in $t$ in the form $(V\,r)$ (resp $(r\,\varepsilon)$). The
 extended structural reduction of $t$ along ${\cal O}$ is the
 transformation to a term $t'$ obtained from $t$ by replacing each
 indexed term $v_{_{\cal O}}$ (the acceptors of ${\cal O}$) by
 $(V\,v)$ (resp $(v\,\varepsilon)$) and erasing the occurence of  $V$
 (resp $\varepsilon$)  in  $(V\,r)$ (resp $(r\,\varepsilon)$) 
. This reduction is denoted by $t \p_{\cal O} t'$.
\end{df}

\begin{re} 
By the definition above, every structural reduction is an
extended structural reduction. It corresponds to the particular
case where the segment-tree consists only of its root.
\end{re}

\begin{ex} Here are two examples of segment-trees and the extended structural reduction. Let $t=(u\;[x. \m
    a.(a\,\<x, (a\;w)\>), y. v])$ and $V$ a value. 
\begin{enumerate}
\item The set ${\cal O}_1 =\{t\}$ is a segment-tree from $t$ in $t$
itself. The acceptors of ${\cal O}_1$ are $\m a.(a\,\<x, (a\;w)\>)$ and
  $v$. Then $t$ is represented  as follows:

 $t=(u\;[x. (\m a.(a\,\<x, (a\;w)\>))_{{\cal O}_1}, y. v_{{\cal O}_1}])$,
 \\and $(V\;t) \p_{{\cal O}_1} (u\;[x. (V\;\m a.(a\,\<x, (a\;w)\>)), y. (V\;v)])$.

\item The set ${\cal O}_2 =\{t, \m a.(a\,\<x, (a\;w)\>)\}$ is also a
  segment-tree from $t$ in $t$. The acceptors of  ${\cal O}_2$ are
  $\<x, (a\;w)\>$,  $w$
  and $v$. Then $t$ is represented as follows:

 $t=(u\;[x. \m a.(a\;\<x, (a\; w_{{\cal O}_2})\>_{{\cal O}_2}),
 y. v_{{\cal O}_2}])$,\\ and $(V\;t) \p_{{\cal O}_2}(u\;[x. \m
 a.(a\;(V\;\<x, (a\;(V\;w))\>)), y. (V\;v)])$. 
\end{enumerate}

\end{ex}

\begin{df}\label{par} 
The parallel reduction $\p$ is defined inductively by the following
rules: 
\begin{itemize}

\item $x \p x$

\item If $t \p t'$, then $\l x.t \p \l x.t'$, $\m a.t \p \m a.t'$,
  $(a\;t) \p (a\;t')$ and $\omega_i t \p \omega_i t'$
 
\item If $t \p t'$ and $u \p u'$, then $\<t,u\> \p \<t',u'\>$

\item If $t \p t'$ and $\varepsilon \tilde{\p} \varepsilon'$,
then $(t \; \varepsilon) \p (t' \; \varepsilon')$

\item If $t \p t'$ and $V \p V'$, then $(\l x.t\; V) \p t'[x:=V']$

\item If $V_i \p V'_i$, then $(\<V_1,V_2\>\; \pi_i) \p V'_i$

\item If $ V \p V'  $ and $u_i \p u_i'$, then $(\omega_iV
  \;[x_1,u_1,x_2,u_2]) \p u_i'[x_i:=V']$ 

\item If  $t \p t'$, $ V \p V'$  (resp $\varepsilon \tilde{\p}
  \varepsilon'$), and ${\cal O}$ is a segment-tree from $t$ in $t$,
  and  $(V'\,t') \p_{\cal O} w$ (resp $(t'\,\varepsilon') \p_{\cal O}
  w$), then $(V\,t) \p w$  (resp $(t\,\varepsilon) \p w$), where $\varepsilon \tilde{\p} \varepsilon'$ means that:

\begin{itemize}
\item $\varepsilon = \varepsilon'= \pi_i$, or
\item ($\varepsilon = u $ and $ \varepsilon' =u'$) or ($\varepsilon
  =[x.u,y.v] $ and $ \varepsilon'=[x.u',y.v']$) such that $u\p u'$ and
  $v \p v'$.
\end{itemize}

\end{itemize}

It is easy to see that $\red_v$ is the transitive closure of $\p$.
\end{df}

\begin{df}\label{devlp}
Let $t$ be a term , we define the complete development $t^*$ as
follows:

\begin{itemize}

\item $x^*=x$
\item $(\l x.t)^*=\l x.t^*$
\item $(\m a.t)^*=\m a.t^*$
\item $\<t_1,t_2\>^* =\<t_1^*,t_2^*\>$
\item $(\o_i t)^*=\o_i t^*$
\item $(a\,t)^*=(a\, t^*)$,
\item $(t\; \varepsilon)^*=(t^*\, \varepsilon^*)$, if $(t\, \varepsilon)$ is not a redex
\item $(\l x.t\; V)^*=t^*[x:= V^*]$
\item $(\<V_1,V_2\> \; \pi_i)^*= V_i^*$
\item $(\o_i V\;[x_1.u_1,x_2.u_2])^*=u_i^*[x:=V^*]$ 
\item Let ${\cal O}_{m}$ be the maximal segment-tree from $t$ in $t$,
  and   $(V^*\;t^*) \p_{{\cal O}_{m}} w$ (resp $(t^*\;\tilde{\varepsilon^*}) \p_{{\cal O}_{m}}
  w$), then $(V\;t)^*= w$ (resp $ (t\;
  \varepsilon)^* =w$), where $\tilde{\varepsilon^*}$ means:
\begin{itemize}
\item $\varepsilon$, if $\varepsilon= \pi_i$
\item $u^*$, if $\varepsilon= u$
\item $[x.u^*,y.v^*]$, if $\varepsilon=[x.u,y.v]$
\end{itemize}
\end{itemize}
\end{df}

\begin{lm}
\begin{enumerate}
\item If $t \p t'$ and $V \p V'$, then $t[x:=V] \p t'[x:=V']$.
\item If $t \p t'$ and $\varepsilon \p \varepsilon'$, then
  $t[a:=^*\varepsilon] \p t'[a:=^*\varepsilon']$.
\item  If $t \p t'$ and $V \p V'$, then
  $t[a:=_*V] \p t'[a:=_*V']$.
\end{enumerate}

\end{lm}

\begin{proof} By a straightforward induction on the structure of $t \p
  t'$.\end{proof}

\begin{lm}(The key lemma) \label{key} If $t \p  t'$, then $ t' \p t^*$.
\end{lm}
 
\begin{proof} The proof of this lemma will be the subject of the
next section.\end{proof}

\begin{theo}(The Diamond Property)
If $t \p t_1$ and $t  \p t_2$, then there exists $t_3$ such that $t_1
\p t_3$ \\ and $t_2  \p t_3$.
\end{theo}

\begin{proof}
It is enough to take $t_3 = t^*$, then theorem holds by the key lemma.
\end{proof}

Since $\red_v$ is identical to the transitive closure of $\p$, we have
the confluence of the call-by-value \\$\l \m^{\wedge \vee} $-calculus.

\begin{theo} 
If $t \red_v t_1$ and $t \red_v t_2$, then there exists a
term $t_3$ such that  $t_1 \red_v t_3$ and $t_2  \red_v t_3$.
\end{theo}

\section{Proof of the key lemma}

For technical reasons (see the example below), we start this section by extending the notion of the segment-tree.

\begin{df}
\begin{enumerate}
\item Let $v$ be a subterm in a term $t$, $v$ is called a bud in $t$
  iff $v$ is $t$ itself or $v$ occurs in $t$ in the form $(a\; v)$
  where $a$ is a free variable in $t$.

\item Let ${\cal O}_1$,...,${\cal O}_n$ be
  segment-trees from respectively $r_1,...,r_n$ in a term $t$, and
  ${\cal P}$ a set of buds (possibly empty) in $t$. Then a
  segment-wood  is a pair $\< {\cal O}_1 \cup...\cup{\cal
  O}_n, \cal P\>$ such that:

\begin{itemize}
\item $r_i$ is a bud in $t$ for each $i$,
\item ${\cal O}_1,...,{\cal O}_n$ and ${\cal P}$
  are mutually disjoints.
\end{itemize}  

\item Let ${\cal Q}=\< {\cal O}_1 \cup...\cup{\cal
  O}_n, \cal P\>$ be a segment-wood in $t$, the elements of ${\cal
  O}_1 \cup...\cup{\cal O}_n$ are called trunk-pieces of $\cal Q$,
  and those of $\cal P$ are called proper-buds of $\cal Q$.
\begin{enumerate}
\item We denote by $Bud(\cal Q)$ the set of buds ${\cal P} \cup
  \{r_1,...,r_n\}$ in $t$.

\item An acceptor of a segment-wood $\cal Q$ is either an acceptor of
  ${\cal O}_i$ for some $i$, either a proper-bud. 

\item The acceptors of $\cal Q$ are indexed by $\cal Q$.

\item If the root $r$ of a segment-tree  $\cal O$ in $t$ is a bud in  $t$, then we identify $\cal O$ with the segment-wood $\<{\cal O}, \emptyset\>$.
\end{enumerate}

\item Let ${\cal Q}$ be a segment-wood in $t$, and $s$ a subterm in
  $t$. The restriction of indexed subterms by ${\cal Q}$ to $s$
  constrcuts a segment-wood  in $s$, which we will denote also by
  ${\cal Q}$ if there is no ambiguity.

\end{enumerate}

\end{df}

\begin{re}
\begin{enumerate}

\item If $v$ is a bud in $t$, then $v$ has no segment-predecessor in
  $t$. Therefore any segment-successor is not a bud.

\item Let ${\cal Q}=\< {\cal O}_1 \cup...\cup{\cal
  O}_n, \cal P\>$ be a segment-wood, since a segment-successor is not a bud, then any acceptor of
  any ${\cal O}_i$ is not in $Bud({\cal Q})$. 

\item The two conditions in $(2)$ of the above definition are
  equivalent to the fact that
  all the elements of $\cal P$ and the buds $r_1,...,r_n$ are distincts.

\item If ${\cal O}$ is a segment-tree from $t$ in $t$, and $s$ is a
  subterm in $t$, then the restriction of ${\cal O}$ to $s$ constructs
  a segment-wood in $s$.

\item Proper-buds and trunk-pieces cannot be treated in a uniform
  way, since in a term, what will be indexed are the proper-buds
  themselves and the acceptors of the trunk-pieces, thing which is
  allowed by a formulation which makes difference between these two
  notions.

\end{enumerate}

\end{re}

\begin{df} 
Let $t,\varepsilon$ be $\cal E$-terms, $V$ a value and $\cal Q$ a
segment-wood in $t$, we define the term $t[V/\cal Q]$ (resp
$t[\varepsilon/\cal Q]$) which is obtained from $t$ by replacing each
indexed term $v_{\cal Q}$ (the acceptors of ${\cal Q}$) in $t$ by $(V\,v)$
(resp $(v\,\varepsilon)$).
\end{df}

\begin{re}
It's clear that if  $(V\,t)
\p_{\cal O} w$ (resp $(t \, \varepsilon) \p_{\cal O} w$), then $w=
t[V/\cal O]$ (resp $w= t[\varepsilon/\cal
O]$)  .
\end{re}

\begin{ex} Let $t = \m a.(a\; \m b.(b\;\o_2 \l s.(a\; \o_1 s)))$
  be a term and $r$ the subterm $\m b.(b\;\o_2 \l s.(a\; \o_1 s))$ in
  $t$. We define two segment-trees from $t$ in $t$, ${\cal
  O}_1=\{t\}$  and ${\cal O}_2=\{t,r\}$, observe that the acceptors of ${\cal
  O}_1$ are $r$ and $\o_1 s$, however those of ${\cal
  O}_2$ are $\o_2 \l s.(a\;\o_1 s)$ and $\o_1 s$. The restriction  ${\cal Q}_1$
  (resp ${\cal Q}_2$) of ${\cal O}_1$ (resp ${\cal O}_2$) to $r$  is the
  following segment-wood: ${\cal Q}_1 =\<\emptyset,\{r, \o_1 s\}\>$
  (resp ${\cal Q}_2=\<\{r\},\{\o_1 s\}\>$). Remark also that
  $Bud({\cal Q}_1) = Bud({\cal Q}_2)$ and the set of trunk-pieces of
  ${\cal Q}_1$ is a subset of that of ${\cal Q}_2$. Suppose that $V$ is a value then:
\begin{itemize}
\item $t[V/{\cal Q}_1]= \m a.(a\;(V\; \m b.(b\;\o_2 \l s.(a\; (V\;\o_1 s)))))$.
\item  $t[V/{\cal Q}_2]= \m a.(a\; \m b.(b\;(V\;\o_2 \l s.(a\;
  (V\;\o_1 s)))))$.
\item  $t[V/{\cal Q}_1] \p t[V/{\cal Q}_2]$
\end{itemize}
\end{ex}

\begin{lm}\label{wood} 
Let ${\cal Q}_1$ and ${\cal Q}_2$ be two segment-woods in a term $t$
such that: $Bud({\cal Q}_1)=Bud({\cal Q}_2)$ and the set of all
trunke-pieces of ${\cal Q}_1$ is a subset of that of ${\cal
Q}_2$. Suppose also that $t\p t'$ and $V \p V'$ (resp $\varepsilon \tilde{\p}
\varepsilon'$), then $t[V/{\cal Q}_1] \p
t'[V'/{\cal Q}_2]$  (resp $t[\varepsilon/{\cal Q}_1] \p
t'[\varepsilon'/{\cal Q}_2]$).
\end{lm}

\begin{proof} 
By induction on $t$. We look at the last rule used for $t
\p t'$. We examine only one case. The others are either
treated similarly, either by a straightforward induction.

 $t=(W\;u)$ and $t'=u'[W'/\cal O]$, where
  $\cal O$ is a segment-tree from $u$ in $u$, $u\p u'$ and
  $ W \p W'$.
\begin{itemize}
\item If $t$ is not an acceptor of ${\cal Q}_1$ and then
  nor of ${\cal Q}_2$: By the induction hypothesis, $u[V/{\cal
  Q}_1] \p u'[V'/{\cal Q}_2]$ and $W[V/{\cal
  Q}_1] \p W'[V'/{\cal Q}_2]$. Since  $\cal O$ is a
  segment-tree from $u$ in $u$, we have:\\ $t[V/{\cal
  Q}_1]=(W[V/{\cal Q}_1]\;u[V/{\cal Q}_1]) \p
  u'[V'/{\cal Q}_2][W'[V'/{\cal Q}_2]/{\cal
  O}]$$=u'[W'/{\cal O}][V'/{\cal Q}_2]=t'[V'/{\cal
  Q}_2]$.

\item If $t$ is an acceptor of  ${\cal Q}_1$ but not of ${\cal Q}_2$:
  Let  ${\cal Q}_2=\<{\cal O}_t \cup {\cal O}_{r_1} \cup ...\cup {\cal
  O}_{r_n},\cal P\>$ and ${\cal Q}_2^-=\< {\cal O}_{r_1} \cup ...\cup {\cal
  O}_{r_n},\cal P\>$, where ${\cal O}_s$ denotes a segment-tree from the
  bud $s$ in $t$. By the induction hypothesis, $u[V/{\cal
  Q}_1] \p u'[V'/{\cal Q}_2^-]$ and $W[V/{\cal
  Q}_1] \p  W'[V'/{\cal Q}_2^-]$. Moreover
  $( W'[V'/{\cal Q}_2^-]\;u'[V'/{\cal Q}_2^-])
  \p_{{\cal O}}$\\$ u'[V'/{\cal Q}_2^-][W'[V'/{\cal
  Q}_2^-]/{\cal O}]$. Hence  $(W[V/{\cal
  Q}_1]\;u[V/{\cal  Q}_1]) \p u'[V'/{\cal Q}_2^-][W'[V'/{\cal
  Q}_2^-]/{\cal O}]$.\\ Therefore, $t[V/{\cal
  Q}_1]=(V\;(W[V/{\cal Q}_1] \;u[V/{\cal Q}_1])) \p u'[V'/{\cal Q}_2][W'[V'/{\cal
  Q}_2^-]/{\cal O}][V'/{\cal O}_t]=$\\$ u'[W'/{\cal
  O}][V'/{\cal Q}_2^-][V'/{\cal O}_t]=u'[W'/{\cal
  O}][V'/{\cal Q}_2]$$= t'[V'/{\cal Q}_2]$.

\item If $t$ is an acceptor of ${\cal Q}_1$ and ${\cal Q}_2$, then
  $t[V/{\cal Q}_1]=(V\;(W [V/{\cal Q}_1] \;u[V/{\cal Q}_1])) \p $\\$
  (V'\;u'[V'/{\cal Q}_2][W'[V'/{\cal Q}_2]/{\cal
  O}])= (V'\; u'[W'/{\cal O}][V'/{\cal Q}_2]
  )=t'[V'/{\cal Q}_2]$.
\end{itemize}
\end{proof}

\begin{proof}[of the key lemma] 

By induction on $t$. We look at the last rule used in $t \p t'$. Only
one case is mentioned:
 $t=(V\,u)$ and  $t'=u'[V'/\cal O]$ where $\cal O$ is a
  segment-tree from $u$ in $u$, $u\p u'$ and $ V \p V'$. In this case $t^*=
  u^*[V^*/{\cal O}_m] $, where ${\cal O}_m$  is the maximal
  segment-tree from $u$ in $u$. Therefore, by the previous lemma (it' s clear that ${\cal O}$ and ${\cal O}_m$ as segment-woods
  satisfy the hypothesis of this lemma \ref{wood}) and  the  induction hypothesis, $u'[V'/{\cal O}] \p u^*[V^*/{\cal O}_m]$.

\end{proof}

\section{Future work}
The strong normalization of this system  cannot be
directly deduced from that of $\l \m ^{ \wedge \vee}$-calculus, since
we consider the symmetric structural reductions $\m'_v$ and
$\d'_v$. Even if the strong normalization of $\l \m \m'$-calculus is
well known (see \cite{davnou2}), the presence of
$\m'_v$ and $\d'_v$ complicates the management of the duplication and
the creation of redexes when the other reductions are
considered.

\acknowledgements
\label{sec:ack}
We wish to thank P. De Groote for helpful discussions.

\nocite{*}
\bibliographystyle{abbrvnat}
\bibliography{cbv}
\label{sec:biblio}

\end{document}